\documentclass[10pt,reqno]{amsart}
\usepackage{amssymb}
\def\noi{\noindent}

\newtheorem{Thm}{Theorem}[section]
\newtheorem{Def}[Thm]{Definition}
\newtheorem{Lm}[Thm]{Lemma}

\newtheorem{state}{Theorem}

\def\cal{\mathcal}
\def\Bbb{\mathbb}
\def\mf{\mathfrak}

\def\<{\langle}
\def\>{\rangle}

\def\a{\alpha}
\def\b{\beta}
\def\d{\delta}
\def\D{\Delta}
\def\th{\theta}

\def\l{\lambda}

\def\Re{\Bbb R}
\def\F{\Bbb F}
\def\C{\Bbb C}
\def\Z{\Bbb Z}

\def\Q{\Bbb Q}

\def\R{\widetilde R}

\def\H{\cal H}

\def\I{\cal I}

\def\h{\mf h}
\def\g{\mf g}
\def\hght{{\rm ht}}

\def\W{\widetilde W}

\begin{document}
\title[]
{The Cherednik kernel and generalized exponents}
\author{Bogdan Ion}
\thanks{Department of Mathematics, University of Michigan, Ann
Arbor, MI 48109.}
\thanks{E-mail address: {\tt bogdion@umich.edu}}

\begin{abstract}
We show how the knowledge of the Fourier coefficients of the
Cherednik kernel leads to combinatorial formulas for generalized
exponents. We recover known formulas for generalized exponents
of irreducible representations parameterized by dominant roots, and
obtain new formulas for the generalized exponents for
irreducible representations parameterized by the dominant elements
of the root lattice which are sums of two orthogonal short roots.
\end{abstract}
\maketitle

\thispagestyle{empty}
\section*{Introduction}
Let $\g$ be a complex simple Lie algebra of rank $n$ and denote by
$G$ its adjoint group. The algebra ${S}(\g)$ of complex valued
polynomial functions on $\g$ becomes a graded representation for
$G$. It is known from the work of Kostant \cite{kostant} that if
$\I$ denotes the subring of $G$--invariant polynomials on $\g$
then $S(\g)$ is free as an $\I$--module and is generated by $\H$, the
space of $G$--harmonic polynomials on $\g$ (the polynomials
annihilated by all $G$--invariant differential operators with
constant complex coefficients and no constant term), or
equivalently $S(\g)=\I\otimes \H$. The space of harmonic
polynomials thus becomes a graded, locally finite representation
of $G$; it can equivalently be thought of as the ring of regular
functions on the cone of nilpotent elements in $\g$. If we denote
by $\H^i$ its $i$--th graded piece, and by $V_\l$ the irreducible
representation of $G$ with highest weight $\l$ we can consider the
graded multiplicity of $V_\l$ in $\H$
$$
E(V_\l):=\sum_{0\leq i} {\rm dim_\C}\left({\rm
Hom}_G(V_\l,\H^i)\right)t^i
$$
As a polynomial with positive integer coefficients $E(V_\l)$ can
be written in the form $$E(V_\l)=\sum_{i=1}^{v_\l} t^{e_i(\l)}$$
such that $e_1(\l)\leq e_2(\l)\leq \cdots \leq e_{v_\l}(\l)$ and
$v_\l$ is the multiplicity of the $0$--th weight space of $V_\l$.
The positive integers $e_i(\l)$ were called by Kostant the
generalized exponents of $V_\l$. The terminology is justified by
the fact that the classical exponents of $G$, the numbers $e_1\leq
\cdots \leq e_n$ which appear in the factorization of the
Poincar\' e polynomial of $G$
\begin{equation*}\label{betti}
p_G(t)=\prod_{i=1}^n (1+t^{2e_i+1})
\end{equation*}
coincide with the generalized exponents of the adjoint
representation of $G$.

To further motivate the importance of generalized exponents note
that by \cite{lusztig} and \cite{hesse} the polynomials $E(V_\l)$
are particular examples of Kazhdan--Lusztig polynomials (for the
affine Weyl group associated to the Weyl group of $G$) and
therefore of considerable combinatorial complexity. The results of
Lusztig and Hesselink describe $E(V_\l)$ as a $t$--analogue of
the $0$-th weight multiplicity of $V_\l$ via a deformation of
Kostant's weight multiplicity formula introduced by Lusztig.

The problem  of computing the classical exponents of $G$ was
initially motivated by the problem of computing the Betti numbers
of $G$ . It turns out that the classical exponents  admit another
description quite different from the one alluded to above.
It was observed independently by A. Shapiro (unpublished) and R.
Steinberg \cite{steinberg} that if we denote by $h(k)$ the number
of positive roots of height $k$ in the root system associated to
$G$ then the number of times $k$ occurs as an exponent of $G$ is
$h(k)-h(k+1)$. This very simple procedure for computing the
classical exponents was justified by Coleman \cite{coleman} modulo
the empirically observed fact that $2N=nh$ ($N$ is the number of
reflexions in the Weyl group  of $G$ and $h$ is the order of a
special element of the Weyl group called the Coxeter
transformation) and by Kostant \cite{tds} who gave a uniform
proof by studying the decomposition of $\g$ into submodules for
the action of a principal three dimensional subalgebra of $\g$.
There is also a proof of this fact  directly  from
Macdonald's factorization of the Poincar\' e polynomial of the
Weyl group of $G$ \cite{maccoxeter} \cite[Section
3.20]{humphreys}.

The main goal of this paper is to explain how the above
description of the classical exponents and similar descriptions of
generalized exponents can be obtained by analyzing the Fourier
coefficients of  the Cherednik kernel, a certain continuous
function on a maximal torus of $G$. Besides recovering the
formulas for generalized exponents of irreducible representations
parameterized by dominant roots, our main result, Theorem
\ref{finalformula}, describes combinatorially the generalized
exponents for irreducible representations parameterized by
dominant elements $\l$ of the root lattice of $\g$ which are sums
of two orthogonal short roots.

To describe this result we need the following notation.
Let $\l$ be a dominant element of the root lattice of $\g$ which can
be written as a sum of two orthogonal short roots and it is not itself 
a root. For any $\gamma$ in the same Weyl group orbit as $\l$ let $n(\gamma)$
be  the number of (unordered) pairs of positive short orthogonal roots 
which sum up to $\gamma$. Let
$h^{}_\l(k):= h^{\prime}_\l(k)- h^{\prime\prime}_\l(k)$, where
$h^\prime_\l(k)$ is the number of weights of $V_\l$ which have
height $k$ and $h^{\prime\prime}_\l(k)$ is the number of weights
$\gamma$ of $V_\l$  in the same Weyl group orbit as $\l$ and whose
height is $k+n(\gamma)$.

\begin{state}
Let $\l$ be a dominant element of the root lattice of $\g$ which can
be written as a sum of two orthogonal short roots and it is not itself 
a root. With
this notation above, the multiplicity of $V_\l$ in $\H^k$ equals
 $h(k)-h(k+1)$.
\end{state}

Our result suggests that similar formulas for generalized exponents
for other classes of irreducible representations of $G$ are also possible
if one explicitly describes the Fourier coefficients of the Cherednik kernel
parametrized by all weights of the irreducible representation under
consideration. A general technique of inductively computing
the Fourier coefficients of the Cherednik kernel is described
in Theorem \ref{prop1.15}. Another closely related method for computing
the Fourier coefficients of the Cherednik kernel was introduced by
Bazlov \cite{bazlov}. It is  based on Cherednik
operators and was succesfuly applied to
compute the Fourier coefficients parametrized by roots, but this
method seems to be less efficient in general because of  the
complexity of Cherednik operators.

\section{Preliminaries}
\subsection{}

Let $\g$ be a complex simple Lie algebra of rank $n$ and denote by
$G$ its adjoint group. Let $\h$ and $\mf b$ be  a Cartan
subalgebra respectively a Borel subalgebra of $\g$ such that
$\h\subset {\mf b}$, fixed once and for all. The maximal torus of
$G$ corresponding to $\h$ is denoted by $H$. We have $H=TA$ where
$T$ is a compact torus and $A$ is a real split torus. The volume
one Haar measure on $T$ is denoted by $ds$.

Let $R\subset \h^*$ be the set of roots of $\g$ with respect to
$\h$, let $R^+$ be the set of roots of $\mf b$ with respect to
$\h$ and denote by $R^-=-R^+$. Of course, $R=R^+\cup R^-$; the
roots in $R^+$ are called positive and those in $R^-$ negative.
The set of positive simple roots determined by $R^+$ is denoted by
$\{\a_1,\dots,\a_n \}$. We know that the roots in $R$ have at most
two distinct lengths. We will use the notation $R_s$ and $R_\ell$
to refer respectively to the short roots and the long roots in
$R$. If the root system is simply laced we consider all the roots
to be short. The dominant element of $R_s$ is denoted by $\th_s$
and the dominant element of $R_\ell$ is denoted by $\th_\ell$.

Any element $\a$ of $R$ can be written uniquely as a sum of simple
roots $\sum_{i=1}^n a_i\a_i$. The height of the root $\a$ is
defined to be $$\hght(\a)=\sum_{i=1}^n a_i.$$ The root of $R$ with
has the largest height is denoted by $\th$. By the above
convention, if $R$ is simply laced then $\th=\th_s$ and if $R$ is
not simply laced then $\th=\th_\ell$.

Denote by $r$ the maximal number of laces in the Dynkin diagram
associated to $\g$. There is a canonical positive definite
bilinear form $(\cdot,\cdot)$ on $\h^*_\Re$ (the real vector space
spanned by the roots) normalized such that $(\a,\a)=2$ for long
roots and $(\a,\a)=2/r$ for short roots.  For any root $\a$ define
$\a^\vee=2\a/(\a,\a)$. We know from the axioms of a root system
that $(\a,\b^\vee)$ is an integer for any roots $\a$ and $\b$. In
fact, the only possible values for $|(\a,\b^\vee)|$ are $0$, $1$
or $2$ if the length of $\a$ does not exceed the length of $\b$
(the value $2$ is attained only if $\a=\pm\b$) and $0$, $r$ if the
length of $\a$ is strictly larger than the length of $\b$.

Define $\rho=\frac{1}{2}\sum_{\a\in R}\a^\vee$. With this notation
the height of any root $\a$ can be written as
$\hght(\a)=(\a,\rho)$. The root lattice $Q$ is the integral span
of the simple roots. For an element $\l$ in $Q$ we define its
height as $\hght(\l)=(\l,\rho)$.


\subsection{}
For any root $\a$ consider the reflexion of the Euclidean space
$\h^*_\Re$ given by $$s_\a(x)=x-(x,\a^\vee)\a .$$ The Weyl group
$W$ of the root system $R$ is the subgroup of ${\rm GL}(\h^*_\Re)$
generated by the reflexions $s_\a$, for all roots $\a$ (the simple
reflexions $s_i:=s_{\a_i}$, $1\leq i\leq n$, are enough). The
scalar product on $\h^*_\Re$ is equivariant with respect to the
action of $W$.

We can extend the bilinear form on $\h^*_\Re$ to a bilinear form
on the real vector space $V:=\h^*_\Re + \Re\d$ by setting
$(\d,V)=0$. The affine root system $\R$ is defined as
$$
\R=\{\a+k\d\  |\  \a\in R,\ k\in \Z  \}.
$$
The set of affine positive roots $\R^+$ consists of affine roots
of the form $\a+k\d$ such that $k$ is positive if $\a$ is a
positive root, and $k$ strictly positive if $\a$ is a negative
root. The affine simple roots are $\a_i$ ($1\leq i\leq n$), and
$\a_0:=\d-\th$.

The affine Weyl group $\W$ is the subgroup of ${\rm GL}(V)$
generated by all reflexions $s_{\a+k\d}$ associated to affine
roots. As above, the affine Weyl group is generated by the simple
reflexions $s_0:=s_{\a_0}$, $s_i$ ($1\leq i\leq n$). Let us
describe explicitly the action of the simple affine reflexion
$$
s_0(x)=s_\th(x)+(x,\th)\d .
$$
The bilinear form on $V$ is equivariant with respect to the affine
Weyl group action.

\section{The Cherednik kernel}

\subsection{}

For an element $\l$ of the root lattice  we denote by $e^\l$ the
corresponding character of the compact torus $T$. The trivial
character $e^0$ will be also denoted by $1$. Let $\Z[Q]$ be the
$\Z$--algebra spanned by all such elements (the group algebra of
the lattice $Q$). Note that the multiplication is given by
$e^\l\cdot e^\mu=e^{\l+\mu}$. There is an involution of $\Z[Q]$
given by $\overline{e^\l}=e^{-\l}$. If we set $e^\d=q$, for $q$ a
fixed complex number, the affine Weyl group acts naturally on
$\C[Q]$. For example, $e^{s_0(\l)}=q^{(\l,\th)}e^{s_\th(\l)}$.

The subalgebra of $\Z[Q]$ consisting of $W$--invariant elements is
denoted by $\Z[Q]^W$.  The irreducible finite dimensional
representations of $G$ are parameterized by the dominant elements
of the root lattice. For a dominant $\l$ we denote by $\chi_\l$
the character of the corresponding irreducible representation of
$G$. Restricting the characters to $T$ we will regard them as
elements of $\Z[Q]$. A basis of $\Z[Q]^W$ is then given by the all
the irreducible characters $\chi_\l$ of $G$.

For any continuous function $f$ on the torus $T$, its Fourier
coefficients are parameterized by $Q$ and are given by
$$
f_\l:=\int_T fe^{-\l}ds .
$$
The  coefficient $f_0$ is called the constant term of $f$; it will
be also denoted by $[f]$.
\subsection{} Let us consider the following function on the torus
$$\D=\frac{1}{|W|}\prod_{\a\in R}(1-e^{\a}).$$
Modulo the normalization by the cardinal of the Weyl group (which
makes the constant term $[\D]=1$) this is the square absolute
value of the Weyl denominator of the root system $R$. The scalar
product on $\Z[Q]^W$ given by
$$
\<f,g\>:=\int_T f\overline{g}\D ds
$$
makes the characters $\chi_\l$ orthonormal.

Assume $q$ and $t$ are complex numbers of small absolute value and
let
$$\nabla(q,t)=\prod_{\a\in
{R}}\prod_{i\geq 0}\frac{1-q^ie^{\a}}{1-q^ite^{\a}}$$ Since $q$
and $t$ are small the infinite product is absolutely convergent
and $\nabla(q,t)$ should be seen as a continuous function on the
torus $T$. In the special case when $t=q^k$ and $k$ is a positive
integer this function was introduced by Macdonald in
\cite{macorthogonal} (see also \cite{macbook}) and used to define
a family of orthogonal polynomials associated to root systems and
depending on the parameters $q$ and $t$. Note that in this case
$\nabla(q,q^k)$ is given by a finite product and no convergence
problems appear; therefore $q$ is not required to have small
absolute value and it can be regarded as a parameter. The constant
term of $\nabla(q,t)$ was subject to conjectures of Macdonald,
later to be proved by Cherednik \cite{cherednik}. The function
$$\D(q,t)=\frac{\nabla(q,t)}{[\nabla(q,t)]}$$
is a $W$--invariant continuous function on the torus with constant
term equal to one. It is also invariant under the transformation
which sends $e^\l$, $q$ and $t$ to their inverses and therefore
well defined also for $q$ and $t$ in a neighborhood of infinity.
We can define the following non--degenerate scalar product on
$\Z[Q]^W$
$$
\<f,g\>_{q,t}^\D:=\int_T f\overline{g} \D(q,t)ds
$$
\subsection{}
Let us consider also the continuous function on $T$ given by
$$K(q,t)=\prod_{\a\in {R^+}}\prod_{i\geq0}
\frac{(1-q^ie^{\a})(1-q^{i+1}e^{-\a})}{(1-tq^ie^{\a})(1-tq^{i+1}e^{-\a})}
$$ Note that with the notation $e^\d=q$ the above function can be
written as
$$
K(q,t)=\prod_{\a\in \R^+}\frac{1-e^{\a}}{1-te^{\a}}
$$
For $t=q^k$ and positive integral $k$ this function first appeared
in Cherednik's work \cite{cherednik} on the Macdonald constant
term conjecture. Unlike $\nabla(q,t)$ it is not invariant under
the Weyl group. The function
$$C(q,t)=\frac{K(q,t)}{[K(q,t)]}$$
is a function on the torus with constant term equal to one, which
we will call the Cherednik kernel. The following result
establishes a very important property of $c_\l(q,t)$, the Fourier
coefficients of $C(q,t)$.
\begin{Thm}{\rm(\cite[(5.1.10)]{macbook}).}
With the above notation, the Fourier coefficients $c_\l(q,t)$ of
$C(q,t)$ are rational functions in $q$ and $t$. Furthermore, the
Cherednik kernel is invariant under the transformation which sends
$e^\l$, $q$ and $t$ to their inverses.
\end{Thm}
Consider now $q$ and $t$ as formal variables and define the field
$\F:=\Q(q,t)$. We extend the involution on $\Z[Q]$ to the group
algebra $\F[Q]$ by setting $\overline{q}=q^{-1}$ and
$\overline{t}=t^{-1}$. Since  $c_\l(q,t)$ are rational functions
in $q$ and $t$ and therefore defined for generic $q$ and $t$ we
can regard them as elements of $\F$. The invariance of the
Cherednik kernel from the above Theorem can be restated as
\begin{equation}\label{eqconjugate}
{c_\l(q^{-1},t^{-1})}=c_{-\l}(q,t)
\end{equation}
We can define the following non--degenerate scalar product on
$\F[Q]$
$$
\<f,g\>_{q,t}^C:=\int_T f\overline{g} C(q,t)ds
$$
For example $c_\l(q,t)=\<1,e^\l\>^C_{q,t}$.  The scalar product
has the property that
$$\<g,f\>^C_{q,t}=\overline{\<f,g\>}^C_{q,t}$$

It is known (see e.g. \cite[(5.1.35)]{macbook}) that the two
scalar product coincide for all elements $f,g\in \F[Q]^W$
\begin{equation}\label{equalscalar}
\<f,g\>_{q,t}^\D=\<f,g\>_{q,t}^C
\end{equation}
We will use the notation $\<\cdot,\cdot\>_{q,t}$ to refer to
$\<\cdot,\cdot\>^C_{q,t}$. It follows from the above relation that
$\D(q,t)$ has also Fourier coefficients which are rational
functions of $q$ and $t$ and therefore regarded as elements of
$\F$.
\subsection{} For each simple affine root consider the following
operator, called Demazure--Lusztig operator, acting on $\F[Q]$ as
follows
$$
T_i(e^\l)=e^{s_i(\l)} + (1-t)\frac{e^\l-e^{s_i(\l)}}{1-e^{-\a_i}}
$$
The following result is due to Cherednik.
\begin{Thm}[\cite{cherednik}] The
Demazure--Lusztig operators are unitary for the above scalar
product on $\F[Q]$. This means that for any $f,g\in \F[Q]$ we have
$$
\<T_i(f),T_i(g)\>_{q,t}=\<f,g\>_{q,t}
$$
\end{Thm}
The unitarity of the Demazure--Lusztig operators will be our main
tool for computing some of the Fourier coefficients of the
Cherednik kernel.

\subsection{} Consider now $t$ as a formal variable. The graded torus
character of $S(\g)$, the algebra of complex valued polynomial
functions on $\g$, is easily seen to be
$${\rm ch}_{S(\g)}(t)=\prod_{\a\in R}\frac{1}{1-te^\a}$$
If $\chi_\l$ denotes the character of the irreducible representation
of $G$ with highest weight $\l$, then the graded multiplicity of
$V_\l$ inside $S(\g)$ can be computed as
\begin{equation}\label{symm-multiplicity}
\<{\rm ch}_{S(\g)}(t),\chi_\l\>=
\frac{1}{|W|}\int_T \nabla(0,t)\overline{\chi}_\l
\end{equation}
As mentioned in  Introduction  if $\I$ denotes the subring of
$G$--invariant polynomials on $\g$
 and
$\H$ the space of $G$--harmonic polynomials on $\g$
then $S(\g)=\I\otimes \H$ as graded $G$--modules. If follows that if
we want to compute $E(V_\l)$,
the graded multiplicity of $V_\l$ inside $\H$, then
we would have to factor out in formula (\ref{symm-multiplicity}) the
graded multiplicity of the trivial representation inside $S(\g)$, or
equivalently the constant term of $\nabla(0,t)$.  We can conclude that
\begin{equation}
E(V_\l)=\<1,\chi_\l \>^\D_{0,t}
\end{equation}
By  formula (\ref{equalscalar}) we can thus express $E_\l$ as a
sum of weight multiplicities of $V_\l$ times values of Fourier
coefficients of the Cherednik kernel at $q=0$. The non--symmetry
of the Cherednik kernel allows various Fourier coefficients
parameterized by elements in the same Weyl group orbit to behave
differently and therefore to contribute differently to the above
scalar product. This feature is not present for the Macdonald
kernel $\D(q,t)$. We will return to the problem of computing the
Fourier coefficients of the Cherednik kernel after some
combinatorial considerations which will allow us to describe them
in simple terms for elements of several Weyl group orbits .

\section{The height function and the Bruhat order}

\subsection{}
For each $w$ in $W$ let $\ell(w)$ be the length of a reduced (i.e.
shortest) decomposition of $w$ in terms of the $s_i$. We have $
\ell(w)=|\Pi(w)| $ where $$ \Pi(w)=\{\a\in R^+\ |\ w(\a)\in
R^-\}.$$ We also denote by $ {}^c\Pi(w)=\{\a\in R^+\ |\ w(\a)\in
R^+\}.$ If $w=s_{j_p}\cdots s_{j_1}$ is a reduced expression of
$w$, then
$$
\Pi(w)=\{\a^{(i)}\ |\ 1\leq i\leq p\},
$$
with $\a^{(i)}=s_{j_1}\cdots s_{j_{i-1}}(\a_{j_i})$.

For each element $\l$ of $Q$ define $\l_+$ to be the unique
dominant element in $W\l$, the orbit of $\l$. Let $w_\l\in W$ be
the unique minimal length element such that $w_\l(\l_+)=\l$.
\begin{Lm}\label{lemma1}
With the notation above, we have
 $$ \Pi(w_\l\hspace{-0.2cm}^{-1})=\{\a\in R^+\ |\ (\l,\a)<0
\}. $$
\end{Lm}
\begin{proof}
Let $\a$ be an element of $\Pi(w_\l\hspace{-0.2cm}^{-1})$. Then
$w_\l^{-1}(\a)$ is a negative root and in consequence
\begin{equation}\label{l221}
0\geq
\left(\l_+,w_\l^{-1}(\a)\right)=\left(w_\l(\l_+),\a\right)=(\l,\a).
\end{equation}
Let us see that above we cannot have equality. If
$w_\l^{-1}=s_{j_p}\cdots s_{j_1}$ is a reduced expression, then
$$
\a\in\Pi(w_\l^{-1})=\{\a^{(i)}\ |\ 1\leq i\leq p\},
$$
with $\a^{(i)}=s_{j_1}\cdots s_{j_{i-1}}(\a_{j_i})$. Suppose that
$$
0=(\l,\a^{(i)})=(\l,s_{j_1}\cdots s_{j_{i-1}}(\a_{j_i}))=
(s_{j_{i-1}}\cdots s_{j_1}(\l),\a_{j_i})
$$
then
$$
s_{j_i}s_{j_{i-1}}\cdots s_{j_1}(\l)=s_{j_{i-1}}\cdots
s_{j_1}(\l),
$$
fact which contradicts the minimality of $w_\l^{-1}$.

Conversely, if the inequality $(\l,\a)<0$ holds for  a positive
root $\a$ then equation (\ref{l221}) shows that $w_\l^{-1}(\a)$ is
a negative root.
\end{proof}


\subsection{} The Bruhat order is a partial order on any Coxeter
group defined in way compatible with the length function. For an
element $w$ we put $w<s_iw$\  if and only if \
$\ell(w)<\ell(s_iw)$. The transitive closure of this relation is
called the Bruhat order. The terminology is motivated by the way
this ordering arises for Weyl groups in connection with inclusions
among closures of Bruhat cells for a corresponding semisimple
algebraic group.

For the basic properties of the Bruhat order we refer to Chapter 5
in \cite{humphreys}. Let us list a few of them (the first two
properties completely characterize the Bruhat order):
\begin{enumerate}
\item For each $\a\in R^+$ we have $s_\a w<w$ if and only if $\a$
is in $\Pi(w^{-1})$ ; \item $w'< w$ if and only if $w'$ can be
obtained by omitting some factors in a fixed reduced decomposition
of $w$ ; \item if $w' \leq w$ then either $s_i w' \leq w$ or $s_i
w' \leq s_iw$ (or both).
\end{enumerate}
We can use the Bruhat order on $W$ do define a partial order on
each orbit of the Weyl group action on $Q$ as follows.
\begin{Def}
Let $\l$ and $\mu$ be two elements of the root lattice which lie
in the same orbit of $W$. By definition $\l<\mu$ if and only if
$w_\l<w_\mu$.
\end{Def}
By the above Definition the dominant element of a $W$--orbit is
the minimal element of that orbit with respect to the Bruhat
order.
\begin{Lm}\label{lemma2}
Let $\l$ be an element of the root lattice such that $s_i(\l)\neq
\l$ for some $1\leq i\leq n$. Then $w_{s_i(\l)}=s_iw_\l$.
\end{Lm}
\begin{proof}
Because $\ell(s_iw_\l)=\ell(w_\l)\pm 1$ and
$\ell(s_iw_{s_i(\l)})=\ell(w_{s_i(\l)})\pm 1$ we have  four
possible situations depending on the choice of the signs in the
above relations. The choice of a plus sign in both relations
translates into $\a_i\not\in\Pi(w_\l^{-1})$ and
$\a_i\not\in\Pi(w_{s_i\cdot\l}^{-1})$ which by Lemma \ref{lemma1}
and our hypothesis implies that $(\a_i,\l)>0$ and $(\a_i,s_i(\l))>
0$ (contradiction). The same argument shows that the choice of a
minus sign in both relations is impossible. Now, we can assume
that $\ell(s_iw_\l)=\ell(w_\l)+ 1$ and
$\ell(s_iw_{s_i(\l)})=\ell(w_{s_i(\l)})- 1$, the other case being
treated similarly. Using the minimal length properties of $w_\l$
and $w_{s_i(\l)}$ we can write
$$
\ell(w_\l)+ 1=\ell(s_iw_\l)\geq \ell(w_{s_i(\l)})=
\ell(s_iw_{s_i(\l)})+1\geq \ell(w_\l)+ 1
$$
which shows that $\ell(s_iw_\l)=\ell(w_{s_i(\l)})$. Our conclusion
now follows from the uniqueness of the element $w_{s_i(\l)}$.
\end{proof}
An immediate consequence is the following
\begin{Lm}\label{lemma1.3}
Let $\l$ be a weight such that $s_i(\l)\neq \l$ for some $1\leq
i\leq n$. Then $s_i(\l)>\l$ if and only if $(\a_i,\l)> 0$. If the
equivalent conditions hold we also have
$$\Pi(w_{s_i(\l)})=\Pi(w_\l)\cup\{w_\l^{-1}(\a_i)\}.$$
\end{Lm}
\begin{Lm}
For an element $\l$ in the root lattice we have
$$
\hght(\l_+)-\hght(\l)=\sum_{\a\in \Pi(w_\l)}(\l_+,\a^\vee)
$$
Moreover, the number $\hght(\l_+)-\hght(\l)-\ell(w_\l)$ is a
positive integer.
\end{Lm}
\begin{proof}
Since $\hght(\l)=(\l,\rho)=(\l_+,w_\l^{-1}(\rho))$ we obtain that
$$\hght(\l_+)-\hght(\l)=(\l_+,\rho-w_\l^{-1}(\rho))$$
If we write
$$
\rho=\frac{1}{2}\sum_{\a\in
\Pi(w_\l^{-1})}\a^\vee+\frac{1}{2}\sum_{\a\in
{}^c\Pi(w_\l^{-1})}\a^\vee
$$
using the equalities
\begin{equation}\label{eq2}
w_\l^{-1}\left(\Pi(w_\l^{-1})\right)=-\Pi(w_\l) \ \text{ and }\
w_\l^{-1}\left({}^c\Pi(w_\l^{-1})\right)={}^c\Pi(w_\l)
\end{equation}
we find that
$$
w_\l^{-1}(\rho)=-\frac{1}{2}\sum_{\a\in
\Pi(w_\l)}\a^\vee+\frac{1}{2}\sum_{\a\in {}^c\Pi(w_\l)}\a^\vee
$$
Our first claim then immediately follows. Regarding the second
claim, note that for $\a\in \Pi(w_\l)$ we always have
$(\l_+,\a^\vee)\geq 1$. Indeed, from the equality (\ref{eq2}) we
know that $\a=-w_\l^{-1}(\b)$ with $\b\in \Pi(w_\l^{-1})$ and
therefore by Lemma \ref{lemma1}
$$
(\l_+,\a^\vee)=-(\l,\b^\vee)>0
$$
In conclusion, $\hght(\l_+)-\hght(\l)-\ell(w_\l)=\sum_{\a\in
\Pi(w_\l)}\left((\l_+,\a^\vee)-1\right)$ is a sum of positive
integers and hence a positive integer.
\end{proof}
For any element $\l$ of the root lattice we will use the notation
$$D_\l=\hght(\l_+)-\hght(\l)-\ell(w_\l)$$
As we will see $D_\l$ encodes a certain type of combinatorial
information about $\l$.

If the root system is not simply laced it will be convenient to
consider
$$D_\l(\ell)=\sum_{\a\in \Pi_\ell(w_\l)}\left((\l_+,\a^\vee)-1\right)$$
and
$$D_\l(s)=\sum_{\a\in \Pi_s(w_\l)}\left((\l_+,\a^\vee)-1\right)$$
where $\Pi_\ell(w_\l)$, respectively $\Pi_s(w_\l)$, is used to
denote the long roots, respectively short roots, in $\Pi(w_\l)$.
\subsection{}
Let us describe $D_\l$ in a few cases. Assume that $\l$ is a short
root. Then $\l_+=\th_s$ and $$D_\l=\sum_{\a\in
\Pi(w_\l)}\left((\th_s,\a^\vee)-1\right)$$ Since $\th_s$ is a
short root, it follows that the scalar product $(\th_s,\a^\vee)$
equals 2 if $\a=\th_s$ and equals 1 otherwise. Therefore $D_\l$
takes the value 1 or 0 depending on whether $\th_s$ is in
$\Pi(w_\l)$ or not. But since $w_\l(\th_s)=\l$ we obtain that
$\th_s$ is in $\Pi(w_\l)$ if and only if $\l$ is a negative root.
Therefore we have proved the following result.
\begin{Lm}\label{lemma1.6}
If $\l$ is a short root then $D_\l=0$ if $\l$ is a positive root
and $D_\l=1$ if $\l$ is a negative root.
\end{Lm}

\subsection{} In the case on non--simply laced root systems we can
investigate $D_\l$ for $\l$ a long root. Denote first by
$N(\th_\ell)$ the number of unordered pairs $\{\a,\b\}$ of short
roots such that $\th_\ell=\a+\b$. For any other long root $\l$ the
number of unordered pairs $\{\a,\b\}$ of short roots such that
$\th_\ell=\a+\b$ is still $N(\th_\ell)$ since $w_\l$ provides a
bijection between the set of such pairs.

If $\a$ and $\b$ are short roots such that $\th_\ell=\a+\b$ then
$(\th_\ell,\a^\vee)=2+(\b,\a^\vee)$. We remark that
$(\th_\ell,\a^\vee)$ cannot be zero and then it equals $r$. It
follows that always $(\b,\a^\vee)=r-2$. The same is true for the
scalar product of pairs of short roots associated in a similar way
to any long root. Denote by $n(\l)$ the number of negative roots
appearing in all unordered pairs of short roots such that
$\l=\a+\b$. The following result describes $D_\l$ in combinatorial
terms.
\begin{Lm}\label{lemma1.7}
For a non--simply laced root system $D_\l(\ell)=0$ if $\l$ is a
positive long root and $D_\l(\ell)=1$ if $\l$ is a negative long
root. Also, $D_\l(s)=(r-1)n(\l)$.
\end{Lm}
\begin{proof}
As before, by examining the scalar products we find that
$D_\l(\ell)=0$ if $\th_\ell$ is not in $\Pi(w_\l)$ and
$D_\l(\ell)=1$ if $\th_\ell$ is in $\Pi(w_\l)$. But since
$w_\l(\th_\ell)=\l$ this translates precisely into our first
claim.

Regarding the second claim we use the fact that
$(\th_\ell,\a^\vee)=r$ for all $\a\in \Pi_s(w_\l)$ to write
$D_\l(s)=(r-1)|\Pi_s(w_\l)|$. Therefore, it will be enough to show
that $n(\l)=|\Pi(w_\l)|$. Remark first that all the unordered
pairs of short root which sum up to $\l$ are of the form
$\{w_\l(\a),w_\l(\b)\}$ with $\a$ and $\b$ positive short roots
such that $\th_\l=\a+\b$. If, for example, $w_\l(\a)$ is a
negative root then $\a\in \Pi_s(w_\l)$. We have shown that
$n(\l)\leq |\Pi(w_\l)|$. For the converse inequality note that if
$\a\in \Pi_s(w_\l)$ then $(\th_\ell,\a)=1$ and hence
$\th_\ell-\a=-s_{\th_\ell}(\a)$ is a short root and
$\th_\ell=\a+(\th_\ell-\a)$. Then $\l=w_\l(\a)+w_\l(\th_\ell-\a)$
and $w_\l(\a)$ is a negative root. In conclusion $|\Pi(w_\l)|\leq
n(\l)$ and our statement is proved.
\end{proof}
\subsection{} We will give a combinatorial description of
$D_\l$ for a few more Weyl group orbits. Let us describe first the
orbits we wish to consider. Define
$$
S:=\{\gamma=\a+\b\ |\ \a,\b\in R_s\ , \ (\a,\b)=0\}
$$
The Weyl group acts on $S$ and the number of orbits of this action
is given by the number of dominant elements of $S$. It is useful
to note that the set $S$ is empty for the root system of type
$G_2$ and that for all the other non--simply laced root systems
the long roots belong to $S$. Let us consider $J$ the set of
connected components of the diagram obtained from the Dynkin
diagram of $R$ by removing the nodes corresponding to those simple
roots for which $(\th_s,\a_i^\vee)=1$ and which contain at least
one node associated to a short simple root of $R$. Note that each
connected component as above is itself a Dynkin diagram and
therefore we can associate its Weyl group $W_j$, root system $R_j$
and highest short root $\th_{s,j}$.
\begin{Lm}
The dominant elements of $S$ are $\th_s+\th_{s,j}$ for all $j$ in
$J$.
\end{Lm}
\begin{proof}

We know (see e.g. \cite{humphreys}) that for any element $x$ of
$\h^*_\Re$ the stabilizer ${\rm stab}_W(x)$ is generated by the
simple reflexions which fix $x$. Therefore the stabilizer of
$\th_s$ is the group generated by the simple reflexions $s_i$ for
which $(\th_s,\a_i)=0$ and using the notation above we obtain that
$${\rm stab}_W(\th_s)=\prod_{j\in J} W_j$$
If for a simple root $\a_i$ we have $(\th_s,\a_i)=0$ then $\a_i$
belongs to one of the root systems $R_j$ and therefore
$(\th_{s,j},\a_i)\geq 0$ for any $j\in J$. Hence,
$(\th_s+\th_{s,j}, \a_i)\geq 0$ for any $j\in J$.  If $\a_i$ is a
simple root such that $(\th_s,\a_i^\vee)=1$, since
$(\th_{s,j},\a_i^\vee)\geq -1$ we obtain again that
$(\th_s+\th_{s,j}, \a_i)\geq 0$ for any $j\in J$. In conclusion,
the elements $\th_s+\th_{s,j}$ are all  dominant. To finish the
proof we will show that any element $\gamma$ of $S$ is in fact
conjugate to one of the $\th_s+\th_{s,j}$.

Fix an element $\gamma=\a+\b$ of $S$ such that $\a,\b\in R_s$ and
$(\a,\b)=0$. We can find a Weyl group element $w$ such that
$w(\a)=\th_s$ and therefore $w(\gamma)=\th_s+w(\b)$. Moreover,
$(\th_s, w(\b))=0$. This means in particular that $s_{w(\b)}$ is
an element of $W$ which fixes $\th_s$. The element $s_{w(\b)}$ of
${\rm stab}_W(\th_s)$ being a reflexion it follows that $w(\b)$ is
a short root in one of the $R_j$. If we denote by $\th_{s,j}$ the
highest short root of $R_j$, we can find an element $w^\prime$ of
$W_j$ such that $w^\prime w(\b)=\th_{s,j}$. Of course, since
$w^\prime$ fixes $\th_s$ we obtain that
$$w^\prime w(\gamma)=\th_s+\th_{s,j}$$ Therefore, we have proved
that each element of $S$ is in the same orbit with one of the
elements $\th_s+\th_{s,j}$, $j\in J$.
\end{proof}
\subsection{}\label{section1.8}
We will investigate the possible values of the scalar products
$(\th_s+\th_s^\prime,\a^\vee)$ for positive roots $\a$. We wish to
study the cases which were not already accounted for. Hence we fix
$j\in J$ such that $\th_s+\th_{s,j}\not = \th_\ell$.

The possible values of the scalar product $(\th_s,\a^\vee)$ are
$0$, $1$ and $2$ and the  possible values of the scalar product
$(\th_{s,j},\a^\vee)$ are $0$, $\pm 1$ and $2$. Note that if one
on the scalar products is $2$ then the other one is necessarily
$0$ since $\a$ is either $\th_s$ or $\th_{s,j}$. Therefore
$(\th_s+\th_{s,j},\a^\vee)=2$ only if $\a=\th_s$, $\a=\th_{s,j}$
or $(\th_s,\a^\vee)=(\th_s^\prime,\a^\vee)=1$. The other possible
values of the scalar product are $0$ and $1$ since
$\th_s+\th_{s,j}$ is dominant and $\a$ positive their scalar
product has to be positive.

The most interesting situation is when we have
$(\th_s,\a^\vee)=(\th_{s,j},\a^\vee)=1$. In the situation when we
have two distinct root lengths  $\a$ can potentially be a long
root. In such a case $(\th_s,\a)=(\th_{s,j},\a)=1$ and therefore
$s_{\th_s} s_{\th_{s,j}}(\a)=\a-r\th_s-r\th_{s,j}$ is a long root.
The scalar product $(r\th_s+r\th_{s,j}-\a,\a)=2(r-1)$. If $r=3$
this leads to a contradiction and if $r=2$ then we obtain that
$\a=\th_s+\th_{s,j}$. Hence, $\a$ being a dominant long root it
must equal $\th_\ell$ and hence $\th_s+\th_{s,j}=\th_\ell$
(contradiction). We have shown that if $\a$ is a positive root and
$(\th_s,\a^\vee)=(\th_{s,j},\a^\vee)=1$ then $\a$ is necessarily
short.

Let $A:=\{\a\in R_s^+\ |\
(\th_s,\a^\vee)=(\th_{s,j},\a^\vee)=1\}$. Denote by
$\varphi=-s_{\th_s}s_{\th_{s,j}}$. We will show that $\varphi$ is
an involution of $A$ without fixed points. Indeed,
$\varphi(\a)=\th_s+\th_{s,j}-\a$ is a short root and
$(\th_s,\th_s+\th_{s,j}-\a)=2/r-1/r=1/r$ and similarly
$(\th_{s,j},\th_s+\th_{s,j}-\a)=1/r$, showing that $\varphi(\a)$
is an element of $A$. Obviously, $\varphi^2$ is the identity. If
$\a$ is fixed by $\varphi$ then $\th_s+\th_{s,j}=2\a$. Computing
the scalar product with $\a$ we obtain $2/r=4/r$ which is a
contradiction. Therefore, the involution $\varphi$ does not have
fixed points.

\subsection{}

One consequence of the above considerations is that $A$ has an
even number of elements. For our $j\in J$ (chosen such that
$\th_s+\th_{s,j}\not =\th_\ell$) denote by $n(j)$ the number of
unordered pairs $\{\a,\b\}$ of short orthogonal roots such that
$\th_s+\th_{s,j}=\a+\b$. Also, for $\l\in W(\th_s+\th_{s,j})$
denote by $n(\l)$ the number of negative roots appearing in all
unordered pairs of short roots such that $\l=\a+\b$.
\begin{Lm}\label{lemma1.9}
With the notation above $n(j)=1+|A|/2$.
\end{Lm}
\begin{proof}
If $\{\a,\b\}$ is  a pair of short orthogonal roots for which
$\th_s+\th_{s,j}=\a+\b$, then
$(\th_s+\th_{s,j},\a^\vee)=(\a+\b,\a^\vee)=2/p$. Such a root must
necessarily be positive since otherwise
$\hght(\a+\b)<\hght(\th_s)$. From previous considerations we know
that either $\a\in A$, either $\a\in\{\th_s,\th_{s,j}\}$.
Therefore, the pair $\{\a,\b\}$ is $\{\th_s,\th_{s,j}\}$ or the
pair $\{\a,\varphi(\a)\}$ for some $\a\in A$. It is easy to see
that the number of such pairs is  $1+|A|/2$.
\end{proof}
The next result describes $D_\l$ in combinatorial terms.
\begin{Lm} For an element $\l\in W(\th_s+\th_{s,j})$ as above we
have $$D_\l=D_\l(s)=n(\l)$$
\end{Lm}
\begin{proof}
As we have argued before, there is no long root $\a$  such that
$(\th_s+\th_{s,j},\a^\vee)=2$ and therefore $D_\l(\ell)=0$.
Furthermore,
$$D_\l=D_\l(s)=\sum_{\a\in
\Pi_s(w_\l)}((\th_s+\th_{s,j},\a^\vee)-1)$$ and since the scalar
product $(\th_s+\th_{s,j},\a^\vee)$ is  at most 2 we obtain that
$D_\l$ is the number of $\a\in \Pi_s(w_\l)$  for which
$(\th_s+\th_{s,j},\a^\vee)=2$. We know from  Lemma \ref{lemma1.9}
that a short positive root $\a$ such that
$(\th_s+\th_{s,j},\a^\vee)=2$ gives rise an expression
$\th_s+\th_{s,j}=\a+\b$ with $\a$ and $\b$ short positive roots
Therefore $\l=w_\l(\a)+w_\l(\b)$ and the short root $w_\l(\a)$ is
negative. We have shown that $n(\l)\geq D_\l$. For the converse
inequality we argue as in the proof of Lemma \ref{lemma1.7}.
\end{proof}
The following result is an immediate consequence of the
combinatorial description of $D_\l$.
\begin{Lm}\label{lemma1.11}
Let $\l\in W(\th_s+\th_{s,j})$ as above such that $\hght(\l)=0$.
Then $D_\l=n(j)$.
\end{Lm}
\begin{proof}
The claim is clear since if $\{\a,\b\}$ is a pair of orthogonal
short roots such that $\l=\a+\b$ then because $\hght(\l)=0$
precisely one of $\a$ or $\b$ is a positive root and the other is
a negative root. In conclusion $n(\l)=n(j)$.
\end{proof}
The next result will be useful later.
\begin{Lm}\label{lemma1.12}
Let $\l=s_\th(\th_s+\th_{s,j})$. Then $D_\l=2n(j)-1$ is the root
system $R$ is simply laced and $D_\l=n(j)$ if the root system is
not simply laced.
\end{Lm}
\begin{proof}
Consider a pair $\{\a,\b\}$ of short orthogonal positive roots
such that $\a+\b=\th_{s}+\th_{s,j}$. Then
$\{s_\th(\a),s_\th(\b)\}$ is pair of short orthogonal roots such
that $s_\th(\a)+s_\th(\b)=\l$ and all the pairs with this property
arise in this way.

Assume first that $R$ is simply laced. If
$\{\a,\b\}=\{\th_s,\th_{s,j}\}$ then
$\{s_\th(\a),s_\th(\b)\}=\{-\th_s,\th_{s,j}\}$. In all the other
cases $\{s_\th(\a),s_\th(\b)\}=\{\a-\th_s,\b-\th_{s}\}$. Since
$\th_s$ is the highest root of $R$ the  number of negative roots
appearing in the $n(j)$ unordered pairs of short roots which sum
up to $\l$ equals $2n(j)-1$.

If $R$ is non--simply laced then $(\th_s+\th_{s,j},\th)=1$ which
forces of course $(\a+\b,\th)=1$. Because $\th$ is dominant both
$(\a,\th)$ and $(\b,\th)$ are positive integers and therefore one
of them equals 0 (say, the first one) and the other equals 1.
Hence $\{s_\th(\a),s_\th(\b)\}=\{\a-\th,\b\}$. In conclusion the
number of negative roots appearing in the $n(j)$ unordered pairs
of short roots which sum up to $\l$ equals $n(j)$.
\end{proof}
\subsection{} For the root system $R$ we denote by $N(R)$ the
number of positive roots in $R$. Similarly we denote by $N(R_s)$
the number of positive short roots in $R$ and we use corresponding
notation for the root systems $R_j$.
\begin{Lm} With the notation above,
there are exactly $N(R_s)N(R_{s,j})/n(j)$ elements in the orbit
$W(\th_s+\th_{s,j})$.
\end{Lm}
\begin{proof}
For a fixed short root there are exactly $2N(R_{s,j})$  short
roots orthogonal to it  and with the sum in the prescribed orbit
(since this is the situation for $\th_s$). Therefore, the total
number of pairs of orthogonal short roots is $4N(R_s)N(R_{s,j})$.
From all these pairs by taking their sum we obtain  each element
of the orbit $W(\th_s+\th_{s,j})$ exactly $2n(j)$ times. In
conclusion the number of elements in the orbit
$W(\th_s+\th_{s,j})$ has the predicted value.
\end{proof}

\section{Fourier coefficients}

In this section we will describe a general inductive procedure for
computing the Fourier coefficients of the Cherednik kernel and then
we apply it to find explicit formulas for the coefficients
corresponding to elements in a few Weyl group orbits.

\subsection{} Let $\l$ be an element of the root lattice and $\a_i$ a
simple root.
If $(\l,\a_i^\vee)=k>0$ then
$$
T_i(e^\l)=e^{s_i(\l)}
+(1-t)(e^\l+\cdots +e^{\l-(k-1)\a_i})
$$
Note that for $1<j<k$ the element $\l-j\a_i$ is a convex
combination of $\l$ and $s_i(\l)$. Indeed
$$
\l-j\a_i=(1-j/k)\l+j/k
s_i(\l)
$$
In consequence they lie in Weyl group orbits strictly closer to
the origin than the elements in $W\l$. The same is true if
$(\l,\th)=k>0$
$$
T_0(e^\l)=tq^{k}e^{s_\th(\l)} +(t-1)(qe^{\l-\th}+\cdots
+q^{k-1}e^{\l-(k-1)\th})
$$
and for $1<j<k$ the element $\l-j\th$ is a convex combination of
$\l$ and $s_\th(\l)$. Using the unitarity of the Demazure--Lusztig
operators we obtain relations between Fourier coefficients of the
Cherednik kernel.

Using the equality $\<T_i(1),T_i(e^\l)\>_{q,t}=\<1,e^\l\>_{q,t}$
we obtain the following relations
\begin{equation}\label{eq4}
tc_{s_i(\l)}(q,t)-c_\l(q,t)=(1-t)\left( c_{\l-\a_i}(q,t)+\cdots
+c_{\l-(k-1)\a_i}(q,t) \right)
\end{equation}
for all $1\leq i\leq n$ such that $s_i(\l)>\l$. Also if
$(\l,\th)=k>0$ we have
\begin{equation}\label{eq5}
tq^{k}c_{\l}(q,t) - c_{s_\th(\l)}(q,t)=(1-t)\left(
q^{k-1}c_{\l-\th}(q,t)+\cdots +qc_{\l-(k-1)\th}(q,t) \right)
\end{equation}
Fix a non--zero dominant element $\l_+\in Q$ and consider the
homogeneous system associated to the above equations. The unknowns
are $x_\l$ for all $\l\in W\l_+$ and the equations
\begin{eqnarray}\label{eq6}
tx_{s_i(\l)}-x_\l&=&0\ \ \   \text{if}\  1\leq i\leq n\
\text{and}\ s_i(\l)>\l \\
tq^{k}x_{\l} - x_{s_\th(\l)} &= &0 \ \ \ \text{if}\  (\l,\th)=k>0
\label{eq7}
\end{eqnarray}

It is easy to see that from equation (\ref{eq6}) we obtain that
$x_\l=t^{-\ell(w_\l)}x_{\l_+}$. We also have $k:=(\l_+,\th)>0$ and
from equation (\ref{eq7}) we get
$$
tq^kx_{\l_+}-t^{-\ell(w_{s_\th(\l)})}x_{\l_+}=0
$$
which implies that $x_{\l_+}=0$ and therefore $x_\l=0$ for all
$\l\in W\l_+$.
\begin{Thm}\label{prop1.15}
The system given by equations (\ref{eq4}) and (\ref{eq5}) and
$c_0(q,t)=1$ has unique solution.
\end{Thm}
\begin{proof}
We know that the system has at least one solution. Fix now a
non--zero dominant $\l_+\in Q$. Since the homogeneous system given
by equations (\ref{eq6}) and (\ref{eq7}) has a unique solution it
follows that the system given by equations (\ref{eq4}) and
(\ref{eq5}) has at most one solution. Combining these two remarks
it follows that for any $\l\in W\l_+$, $c_\l(q,t)$ is uniquely
expressible in terms of $c_\mu(q,t)$'s, where $\mu$ lies in an
orbit of $W$ closer to the origin that $\l_+$. Using induction on
the distance of $\l_+$ to the origin we obtain that our system has
a unique solution.
\end{proof}
\subsection{} We will now apply the inductive procedure described
in the proof of Theorem \ref{prop1.15} to find the Fourier
coefficients of the Cherednik kernel corresponding to a few Weyl
group orbits. The orbit of $W$ on $Q$ closest to the origin is
$W\th_s$. If the root system is simply laced then $\th_s=\th$ the
highest root of $R$. Let $X_{\th_s}$ be the element of $\F$ for
which
\begin{equation}\label{eq8}
c_{\th_s}(q,t)=t^{\hght(\th_s)}X_{\th_s}
\end{equation}
\begin{Thm}\label{thmshort}
For any  short root $\l$ we have
\begin{equation}\label{eq9}
c_\l(q,t)=t^{\hght(\l)+D_\l}X_{\th_s}+(t^{\hght(\l)}-t^{\hght(\l)+D_\l})
\end{equation}
and $X_{\th_s}=(1-t^{-1})/(1-qt^{\hght(\th)})$.
\end{Thm}
\begin{proof}
We will show that the above formula is valid by induction on the
order on the orbit $W\th_s$ induced from the Bruhat order. For the
minimal element of the orbit we have $D_{\th_s}=0$ (by Lemma
\ref{lemma1.6}) and the predicted formula coincides with
(\ref{eq8}).

Assuming that the predicted formula is true for $\l$ we will show
that it is true for $s_i(\l)>\l$. As explained in Lemma
\ref{lemma1.3} this means that $(\l,\a_i^\vee)>0$. In fact the
possible values of the scalar product are $1$ or $2$ (only if
$\l=\a_i$). If $(\l,\a_i^\vee)=1$ then
$\hght(s_i(\l))=\hght(\l)-1$ and
$\ell(w_{s_i(\l)})=\ell(w_{\l})+1$. It follows that
$D_{s_i(\l)}=D_\l$ therefore using equation (\ref{eq4}) we get
\begin{eqnarray*}
c_{s_i(\l)}(q,t)&=&t^{-1}c_\l(q,t) \\
&=&
t^{\hght(\l)+D_\l-1}X_{\th_s}+(t^{\hght(\l)-1}-t^{\hght(\l)+D_\l-1})\\
&=&
t^{\hght(s_i(\l))+D_{s_i(\l)}}X_{\th_s}+(t^{\hght(s_i(\l))}-t^{\hght(s_i(\l))+D_{s_i(\l)}})
\end{eqnarray*}
If $(\l,\a_i^\vee)=2$ then $\l=\a_i$ and $s_i(\l)=\a_i$. In
consequence $D_\l=0$ and  $D_{s_i(\l)}=1$. Again by equation
(\ref{eq4}) we get
\begin{eqnarray*}
c_{s_i(\l)}(q,t)&=&t^{-1}c_\l(q,t) +(t^{-1}-1)\\
&=&
X_{\th_s}+(t^{-1}-1)\\
&=&
t^{\hght(s_i(\l))+D_{s_i(\l)}}X_{\th_s}+(t^{\hght(s_i(\l))}-t^{\hght(s_i(\l))+D_{s_i(\l)}})
\end{eqnarray*}
We have thus shown that the formula (\ref{eq9}) is valid for all
short roots. To show that $X_{\th_s}$ has the predicted value we
use the equation (\ref{eq5}).

Indeed, if the root system $R$ is simply laced then $\th=\th_s$
and the equation (\ref{eq5}) for $\l=\th_s$ becomes
$$
tq^2c_{\th_s}(q,t)=c_{-\th_s}(q,t)+q(1-t)
$$
and replacing $c_{\th_s}(q,t)$ and  $c_{-\th_s}(q,t)$ with their
formulas in terms of $X_{\th_s}$ we obtain the desired result.

If the root system $R$ is non--simply laced then $\th=\th_\ell$
and the equation (\ref{eq5}) for $\l=\th_s$ becomes
$$
tqc_{\th_s}(q,t)=c_{\th_s-\th}(q,t)
$$
and again our result follows.
\end{proof}

\subsection{} We will describe next the Fourier coefficients of
the Cherednik kernel corresponding to long roots in the case of a
non--simply laced root systems.
\begin{Thm}\label{thmlong}
For any long root $\l$ we have
\begin{equation}\label{eq10}
c_\l(q,t)=t^{\hght(\l)+D_\l(\ell)}X_{\th_s} +
(t^{\hght(\l)}-t^{\hght(\l)+D_\l(\ell)})
\end{equation}
\end{Thm}
\begin{proof}
Let us denote by $X_{\th}$ the element of $\F$ defined by
$c_{\th}(q,t)=t^{\hght(\th)}X_{\th} $. We will show inductively
that for any long roots $\l$ we have
\begin{equation}\label{eq11}
c_\l(q,t)=t^{\hght(\l)+D_\l}X_{\th}
+(t^{\hght(\l)+D_\l(\ell)}-t^{\hght(\l)+D_\l})X_{\th_s}+
(t^{\hght(\l)}-t^{\hght(\l)+D_\l(\ell)})
\end{equation}

The formula (\ref{eq11}) is clearly true for $\th$. Assuming that
the predicted formula is true for $\l$ we will show that it is
true for $s_i(\l)>\l$. For such an $\a_i$ the possible values of
the scalar product $(\l,\a_i^\vee)$ are 1 (if $\a_i$ is a long
root), 2 (if $\a_i=\l$) and $r$ (if $\a_i$ is a short root). We
analyze these cases separately.

If $(\l,\a_i^\vee)=1$ then $\hght(s_i(\l))=\hght(\l)-1$ and
$\ell(w_{s_i(\l)})=\ell(w_{\l})+1$. It follows that
$D_{s_i(\l)}=D_\l$ therefore using equation (\ref{eq4}) we get the
desired formula for $c_{s_i(\l)}(q,t)$.

If $(\l,\a_i^\vee)=2$ then $\l=\a_i$ and $s_i(\l)=\a_i$. In
consequence $D_\l(\ell)=0$, $D_{s_i(\l)}(\ell)=1$ and
$D_{s_i(\l)}(s)=D_\l(s)$ . Again equation (\ref{eq4}) gives the
predicted formula for $c_{s_i(\l)}(q,t)$.

If $\a_i$ is a short root and $(\l,\a_i^\vee)=r$ then
$\l-\a_i,\l-(r-1)\a_i$ are short roots (for example
$\l-\a_i=-s_\l(\a_i)$). From Lemma \ref{lemma1.3} we know that
$$
\Pi(w_{s_i(\l)})=\Pi(w_\l)\cup\{w_\l^{-1}(\a_i)\}
$$
and since $\a_i$ is a short root we obtain
$$\Pi_\ell(w_{s_i(\l)})=\Pi_\ell(w_\l)\ \  \text{and}\ \
\Pi_s(w_{s_i(\l)})=\Pi_s(w_\l)\cup\{w_\l^{-1}(\a_i)\}$$ Therefore,
by Lemma \ref{lemma1.7} it follows that
$$D_{s_i(\l)}(\ell)=D_\l(\ell)\ \  \text{and} \ \ D_{s_i(\l)}(s)=D_\l(s)+r-1$$
One immediate consequence of the above equalities is that
\begin{equation*}
\hght(s_i(\l))+D_{s_i(\l)}=\hght(\l)+D_\l-1
\end{equation*}
and that $\l$ and $s_i(\l)$ are either both positive roots or both
negative roots. Hence $$D_{\l-\a_i}=D_{\l-(r-1)\a_i}=D_\l(\ell)$$
The equation (\ref{eq4}) in this case becomes
\begin{equation*}
c_{s_i(\l)}(q,t)=t^{-1}c_\l(q,t)
+(t^{-1}-1)(c_{\l-\a_i}(q,t)+c_{\l-(r-1)\a_i}(q,t))
\end{equation*}
and the induction hypothesis together with the above equalities
implies the formula (\ref{eq11}) for $c_{s_i(\l)}(q,t)$.

We have thus shown that the formula (\ref{eq11}) is valid for all
long roots. To show that (\ref{eq10}) is true it is enough to show
that $X_\th=X_{\th_s}$. To see that this is indeed the case we use
equation (\ref{eq5}) for $\l=\th$
\begin{equation*}
tq^2c_{\th}(q,t)=c_{-\th}(q,t)+q(1-t)
\end{equation*}
or equivalently
\begin{equation*}
(q^2t^{2\hght(\th)}-t^{D_{-\th}(s)})X_\th=(1-t^{D_{-\th}(s)})X_{\th_s}+
(t^{-1}-1)(qt^{\hght(\th)}+1)
\end{equation*}
Replacing $t^{-1}-1$ with $(qt^{\hght(\th)}-1)X_{\th_s}$ in the
above formula immediately gives $X_\th=X_{\th_s}$ and therefore
the desired result.
\end{proof}
Theorem \ref{thmshort} and Theorem \ref{thmlong} were also obtained by
Bazlov \cite[Theorem 3]{bazlov} using a different (but related) procedure for
computing the Fourier coefficients based on the unitarity of Cherednik
operators. Given the considerable complexity of Cherednik operators,
their action is very hard to be analyzed in general.


\subsection{} We will next describe the Fourier
coefficients for the Cherednik kernel for all $\l\in
W(\th_s+\th_{s,j})$, where $j\in J$ is such that
$\th_s+\th_{s,j}\not = \th_\ell$. First, let us define $X_j$ to be
the unique element of $\F$ for which
\begin{equation}\label{eq12}
c_{\th_s +\th_{s,j}}(q,t)=t^{\hght(\th_s+\th_{s,j})}X_jX_{\th_s}
\end{equation}
We will also need the following notation: for $\l$ as above $d_\l$
is defined to be 0 if $\hght(\l)\geq 0$ and to be 1 if
$\hght(\l)<0$. Also,
\begin{eqnarray*}
a_\l(t)&:=&t^{\hght(\l)}(t^{d_\l}+t^{D_\l-n(j)}-t^{D_\l}-t^{d_\l+D_\l-n(j)})\\
b_\l(t)&:=&t^{\hght(\l)}(1+t^{d_\l+D_\l-n(j)}-t^{d_\l}-t^{D_\l-n(j)})
\end{eqnarray*}
\begin{Thm}
For $\l$  as above we have
\begin{equation}\label{eq13}
c_\l(q,t)=t^{\hght(\l)+D_\l}X_jX_{\th_s}+a_\l(t)X_{\th_s} +b_\l(t)
\end{equation}
and $X_j=(1-t^{-n(j)})/(1-qt^{\hght(\th)-n(j)+1})$
\end{Thm}
\begin{proof}
As before, we will show inductively that the formula (\ref{eq13})
holds for any $\l$ in the Weyl group orbit of $\th_s+\th_{s,j}$.
The formula (\ref{eq13}) is easily seen to be true for the
dominant element of the orbit.

Assuming that the predicted formula is true for $\l$ we will show
that it is true for $s_i(\l)>\l$. According to Section
\ref{section1.8} the possible values of the scalar product
$(\l,\a_i^\vee)$ are 1 and 2. In the latter case $\a_i$ is
necessarily a short root.

If $(\l,\a_i^\vee)=1$ then $\hght(s_i(\l))=\hght(\l)-1$ and
$\ell(w_{s_i(\l)})=\ell(w_{\l})+1$. It follows that
$D_{s_i(\l)}=D_\l$ and by using equation (\ref{eq4}) we get the
desired formula for $c_{s_i(\l)}(q,t)$. In fact, special care is
needed if $\hght(\l)=0$ since $d_\l=0$ and $d_{s_i(\l)}=1$, but
Lemma \ref{lemma1.11} assures that everything goes smoothly.

If $\a_i$ is a short root and $(\l,\a_i^\vee)=2$ then $\l-\a_i$ is
also a short root which is orthogonal on $\a_i$ (see Section
\ref{section1.8}). From Lemma \ref{lemma1.3} we know that
$$
\Pi(w_{s_i(\l)})=\Pi(w_\l)\cup\{w_\l^{-1}(\a_i)\}
$$
and therefore
$$D_{s_i(\l)}=D_\l+1$$

The equation (\ref{eq4}) in this case becomes
\begin{equation*}
c_{s_i(\l)}(q,t)=t^{-1}c_\l(q,t) +(t^{-1}-1)c_{\l-\a_i}(q,t)
\end{equation*}
Let us assume first that $\hght(\l)\not = 0$. Since $\l-\a_i$ is a
short root it does not have height zero and moreover
$D_{\l-\a_i}=d_\l=d_{s_i(\l)}$. Using this equality, the above
equation gives
\begin{eqnarray*}
a_{s_i(\l)}(t)&=&t^{\hght(\l)-1}(t^{d_\l}+
t^{D_\l-n(j)}-t^{D_\l}-t^{d_\l+D_\l-n(j)})+(t^{-1}-1)t^{\hght(\l)-1+d_{\l}}
\\
&=&t^{\hght(\l)-1}(t^{d_{\l}-1}+
t^{D_\l-n(j)}-t^{D_\l}-t^{d_\l+D_\l-n(j)})\\
&=&t^{\hght(s_i(\l))}(t^{d_{s_i(\l)}}+
t^{D_{s_i(\l)}-n(j)}-t^{D_{s_i(\l)}}-t^{d_{s_i(\l)}+D_{s_i(\l)}-n(j)})\\
\end{eqnarray*}
and
\begin{eqnarray*}
b_{s_i(\l)}(t)&=&t^{\hght(\l)-1}(1-t^{d_\l})(1-t^{D_\l-n(j)})
+t^{\hght(\l)-1}(1-t^{d_{\l}})(t^{-1}-1)
\\
&=&t^{\hght(\l)-1}(1-t^{d_{\l}})(t^{-1}-t^{D_\l-n(j)})\\
&=&t^{\hght(s_i(\l))}(1-t^{d_{s_i(\l)}})(1-t^{D_{s_i(\l)}-n(j)})\\
\end{eqnarray*}
which give the desired formula for $c_{s_i(\l)}(q,t)$. The case
$\hght(\l)=0$ is handled similarly with the use of Lemma
\ref{lemma1.11}. We have therefore proved the formula (\ref{eq13})
is always valid. We are left with finding the precise value of
$X_j$. As before, we will make use of equation (\ref{eq5}) for
$\th_s+\th_{s,j}$, which takes different forms depending on
whether the root system $R$ is simply laced or not. In both cases
we will use the Lemma \ref{lemma1.12} which describes
$D_{s_{\th}(\th_{s}+\th_{s,j})}$.

If $R$ is non--simply laced the equation (\ref{eq5}) for
$\l=\th_s+\th_{s,j}$ becomes
\begin{equation*}
qtc_{\th_s+\th_{s,j}}(q,t)=c_{\th_s+\th_{s,j}-\th}(q,t)
\end{equation*}
or equivalently
\begin{equation*}
qt^{\hght(\th_s+\th_{s,j})+1}X_jX_{\th_s}=t^{\hght(\th_s+\th_{s,j}-\th)+n(j)}X_jX_{\th_s}
+t^{\hght(\th_s+\th_{s,j}-\th)}(1-t^{n(j)})X_{\th_s}
\end{equation*}
which gives the desired expression for $X_j$.

If $R$ is simply laced the equation (\ref{eq5}) for
$\l=\th_s+\th_{s,j}$ becomes
\begin{equation*}
q^2tc_{\th_s+\th_{s,j}}(q,t)=c_{-\th_s+\th_{s,j}}(q,t)+q(1-t)c_{\th_{s,j}}(q,t)
\end{equation*}
After replacing each Fourier coefficient with its formula given by
(\ref{eq13}) and (\ref{eq9}) and after straightforward
manipulations we find the predicted expression for $X_j$.
\end{proof}

\subsection{} We will collect all the values of the Fourier
coefficients at $q=0$. For $\l$ a positive root (short or long) we
have
\begin{equation}\label{eq14}
c_\l(0,t)=t^{\hght(\l)}-t^{\hght(\l)-1}
\end{equation}
and for $\l$ a negative root we have $c_\l(0,t)=0$.

For $\l$ an element of the orbit of $\th_s+\th_{s,j}(\not =
\th_\ell)$ we get
\begin{equation}\label{eq15}
c_\l(0,t)=(t^{\hght(\l)}-t^{\hght(\l)-1})-(t^{\hght(\l)+D_\l-n(j)}-t^{\hght(\l)+D_\l-n(j)-1})
\end{equation}
for $\l$ of strictly positive height and $c_\l(0,t)=0$ for $\l$ of
height smaller than zero.

Note that the fact that the Fourier coefficients $c_\l(0,t)$ a
zero if $\l$ has not strictly positive height follows directly
from the definition of $K(0,t)$.

For $\l$ a dominant element of the root lattice and a strictly
positive integer $i$ denote by $h^\prime_\l(i)$ the number of
weights (counted with multiplicities) of $V_\l$ the irreducible
representation of $G$ with highest weight $\l$ which have height
$i$. For $\th_s+\th_{s,j}$ denote by
$h^{\prime\prime}_{\th_s+\th_{s,j}}(i)$ the number of elements
$\gamma$ of the orbit of $\th_s+\th_{s,j}$ for which
$\hght(\gamma)+D_\gamma-n(j)=i$. For $\l$ a dominant root let
$h_\l(i):=h^\prime_\l(i)$; for $\th_s+\th_{s,j}$ let
$h_{\th_s+\th_{s,j}}(i)=h^{\prime}_{\th_s+\th_{s,j}}(i)-h^{\prime\prime}_{\th_s+\th_{s,j}}(i)$.
Note that for the above cases $h_\l(1)\geq h_\l(2)\geq \cdots \geq
h_\l(\hght(\l))$ and $h_\l(i)=0$ for $i$ strictly larger than
$\hght(\l)$.
\begin{Thm}\label{finalformula}
The  generalized exponents of $V_\l$ where $\l$ a dominant root or
$\th_s+\th_{s,j}$ are given by the elements of the partition dual
to $\{h_\l(i)\}_{i\geq 1}$. Equivalently, the multiplicity of
$V_\l$ in $\H^i$ is $h_\l(i)-h_\l(i+1)$.
\end{Thm}
\begin{proof}
The result is a direct consequence of the formulas (\ref{eq14})
and (\ref{eq15}). Indeed, if we denote by ${\rm wt}(\l)$ the set
of weights of $V_\l$ and for any $\gamma\in {\rm wt}(\l)$ we
denote by $m_{\l \gamma}$ the weight multiplicity of $\gamma$ in
$V_\l$,  we have
\begin{eqnarray*}
\<1,\chi_\l\>_{0,t}&=& \sum_{\gamma\in {\rm wt}(\l)}
m_{\l\gamma} c_\gamma(0,t)\\
&=& \sum_{1\leq i\leq \hght(\l)} h_\l(i)(t^i-t^{i-1}) +m_{\l 0}\\
&=& \sum_{1\leq i\leq \hght(\l)} (h_\l(i)-h_\l(i+1))t^i
\end{eqnarray*}
Our result is proved.
\end{proof}
Note that for the representation $V_\th$ (the adjoint representation
of $G$) the above Theorem describes
precisely the procedure for computing the classical exponents proved
by Kostant \cite{tds}.
\subsection{} For each type of root system we will describe next the results
of the description of generalized exponents for
$\l=\th_s+\th_{s,j}$, $j\in J$ provided by Theorem
\ref{finalformula}.

\noi $A_n$: $|J|=1$, $n(j)=2$

\noi $ E(\l)=\sum_{1\leq i< \frac{n}{2}}^{} i(t^{2i}+t^{2n-2i}) +
\sum_{1\leq i< \frac{n-2}{2}}^{}
i(t^{2i+1}+t^{2n-2i-1})+\lfloor\frac{n-2}{2}\rfloor t^n$

\noi $B_n$: $|J|=0$.

\noi $C_n$, $n\geq 4$: $|J|=1$, $n(j)=3$

\noi $ E(\l)=\sum_{1\leq i< \frac{n-1}{2}}^{}
i(t^{4i}+t^{4n-4i-4}) + \sum_{1\leq i< \frac{n-2}{2}}^{}
i(t^{4i+2}+t^{4n-4i-6})+\lfloor\frac{n-3}{2}\rfloor t^{2n-2}$

\noi $D_4$: $|J|=3$, $n(j)=3$ for any $j\in J$. Also for all three
$\l=\th_s+\th_{s,j}$ we have

\noi $E(V_\l)=t^2+t^4+t^6$

\noi $D_n$, $n\geq 5$: $|J|=2$. If we denote $J=\{j_1,j_2\}$ we
have $n(j_1)=n-1$, $n(j_2)=3$. With the notation
$\l_1=\th_s+\th_{s,j_1}$ and $\l_2=\th_s+\th_{s,j_2}$ we have

\noi $E(V_{\l_1})=\sum_{i=1}^{n-1}t^{2i}$

\noi For $n$ even we have

\noi $E(V_{\l_2})=\sum_{i=1}^{n/2-2}
\left(\lfloor\frac{i+1}{2}\rfloor (t^{2i} +t^{4n-2i-8})+
(\frac{n}{2}-\lfloor\frac{i+1}{2}\rfloor)(t^{2n-2i-6}
+t^{2n+2i-2})\right)+\\
\hspace*{1.5cm}\frac{n}{2}(t^{2n-6}+t^{2n-4}+t^{2n-2})$

\noi For $n$ odd we have

\noi $E(V_{\l_2})=\sum_{i=1}^{n-3} \lfloor\frac{i+1}{2}\rfloor
(t^{2i} +t^{4n-2i-8})+\sum_{i=\frac{n-3}{2}}^{n-3}(t^{2i+1}
+t^{4n-2i-9})+ \frac{n-1}{2}t^{2n-4}$

\noi $E_6$: $|J|=1$, $n(j)=4$

\noi $E(V_\l)=t^2+t^3+t^4+t^5+2t^6+t^7+2t^8+2t^9
+2t^{10}+t^{11}+2t^{12}+t^{13}+t^{14}+t^{15}+t^{16} $

\noi $E_7$: $|J|=1$, $n(j)=5$

\noi $E(V_\l)=t^2+t^4+2t^6+2t^8
+3t^{10}+3t^{12}+3t^{14}+3t^{16}+3t^{18}+2t^{20}+2t^{22}+t^{24}+t^{26}
$

\noi $E_8$: $|J|=1$, $n(j)=7$

\noi $ E(V_\l)=t^2+t^6+t^8+t^{10}
+2t^{12}+2t^{14}+t^{16}+3t^{18}+2t^{20}+2t^{22}+3t^{24}+2t^{26}+\hspace*{1.3cm}
 2t^{28}+3t^{30}+t^{32}
+2t^{34}+2t^{36}+t^{38}+t^{40}+t^{42}+t^{46} $

\noi $F_4$: $|J|=0$.

\noi $G_2$: $|J|=0$.

It is interesting to note that whenever it is defined $n(j)$
equals some classical exponent for the root system in question. We
also observe the following symmetry (recall that $v_l$ is the
multiplicity of the $0$--th weight space of $V_\l$)
$$
e_i(\l)+e_{v_\l-i}(\l)=\hght(\th_s)+\hght(\th_{s,j})+2
$$

\subsection{} We close by noting that if $e_i=e_i(\th)$ are the
classical exponents, their  symmetry observed by Chevalley (and
proved by Kostant \cite{tds})
\begin{equation}\label{chevalley}
e_i+e_{n-i}=\hght(\th)+1
\end{equation}
can be also explained as follows. By examining $\D(q,t)$ we
observe that the scalar product
$\<1,\chi_\th\>_{t,t}=\<1,\chi_\th\>= 0$. On the other hand, it is
easy to see that
$$
\<1,\chi_\th\>_{q,t}=\frac{(t^{e_1}+\cdots +
t^{e_n})-qt^{\hght(\th)}(t^{-e_1}+\cdots +
t^{-e_n})}{1-qt^{\hght(\th)}}
$$
and therefore our convention $e_1\leq\cdots\leq e_n$ immediately
gives (\ref{chevalley}). For non--simply laced root systems the
same argument gives a similar symmetry for the generalized
exponents for $V_{\th_s}$. If $n_s$ denotes the number of simple
short roots in $R$ then we have
$$
e_i(\th_s)+e_{n_s-i}(\th_s)=\hght(\th)+1
$$
However, this type of considerations do not explain the symmetry
of the generalized exponents observed in the previous paragraph.



\end{document}